\documentclass[conference,letterpaper]{IEEEtran}
\usepackage[left=19.1mm, right=19.1mm, top=18.5mm, bottom=19.9mm]{geometry}
\IEEEoverridecommandlockouts
\usepackage{cite}
\DeclareFontFamily{OT1}{pzc}{}
\DeclareFontShape{OT1}{pzc}{m}{it}{<-> s * [1.10] pzcmi7t}{}
\DeclareMathAlphabet{\mathpzc}{OT1}{pzc}{m}{it}
\usepackage{amsmath,amssymb,amsfonts,amsthm}
\newtheoremstyle{assumptionstyle}
   {0pt}
   {0pt}
   {}
   {11pt}
   {\em}
   {\em}
   {5pt}
   {}
\theoremstyle{assumptionstyle}
\newtheorem{assumption}{Assumption}
\usepackage{algorithm}
\usepackage{algorithmic}

\usepackage{soul} 

\usepackage{cleveref} 

\usepackage{graphicx}
\usepackage{textcomp}
\usepackage{xcolor}
\def\BibTeX{{\rm B\kern-.05em{\sc i\kern-.025em b}\kern-.08em
    T\kern-.1667em\lower.7ex\hbox{E}\kern-.125emX}}

\usepackage{todonotes}
\usepackage{comment}
\newcommand{\myequation}{\begin{equation}}
\newcommand{\myendequation}{\end{equation}}
\let\[\myequation
\let\]\myendequation
\allowdisplaybreaks
\begin{document}
\title{\vspace*{+6.mm}  Learning Model Predictive Control for Connected Autonomous Vehicles
\thanks{The authors are with the University of Virginia. \{hj2bh,cf5eg\}@virginia.edu. This work was partially supported by the NSF under grants CPS-1739333.}
}

\author{\IEEEauthorblockN{Hassan Jafarzadeh and Cody Fleming}
}

\maketitle

\begin{abstract}
A Learning Model Predictive Controller (LMPC) is presented and tailored to platooning and Connected Autonomous Vehicles (CAVs) applications. The proposed controller builds on previous work on nonlinear LMPC, adapting its architecture and extending its capability to (a) handle dynamic environments and (b) account for data-driven decision variables that  derive from an unknown or unknowable function. The paper presents the control design approach, and shows how to recursively construct an outer loop candidate trajectory and an inner iterative LMPC controller that converges to an optimal strategy over both model-driven and data-driven variables. Simulation results show the effectiveness of the proposed control logic.
\end{abstract}


\section{Introduction}\label{sec:intro}

Recent advances in autonomous driving are becoming increasingly ubiquitous, and 
Advanced Driving Assistance Systems (ADAS) have the potential to improve safety and comfort in various driving conditions. 
In addition, the use of wireless communication networks 
could enable autonomous vehicles to reach high performance states that would not otherwise be feasible or safe (closer distances, higher speeds, etc). However, there is currently no methodology for providing provable safety guarantees for such states. 

The main challenge is to capture the interdependence between mobility, wireless, and safety: a vehicle's motion has a profound effect on the wireless channel, which in turn affects the ability to maintain physical safety. Many current efforts in 
the broader context of networked multi-agent systems ~\cite{nowzari2016distributed,dimarogonas2012distributed,wang2016multi,wang2016safety}  assume a stationary channel: future properties of the wireless channel will be similar to those of the past. However, this assumption ignores the effect of mobility on the wireless channel, and so the controller cannot choose a motion plan that would improve safety (and, as a consequence, performance) by improving the wireless channel. In this paper, We explore these notions in the context of so-called {vehicle platooning}. 

Adaptive cruise control (ACC) is a widely used ADAS module that controls the vehicle longitudinal dynamics. ACC is triggered once a preceding vehicle is detected within a certain distance range from the ego vehicle. ACC automatically maintains a proper minimum safe distance from preceding vehicles by automatically adjusting braking and acceleration. ACC enhances mobility, improves safety and comfort, and reduces energy consumption. The use of Model Predictive Control (MPC) for ACC applications is becoming increasingly common in the literature \cite{luo2010model,corona2008adaptive}. 

The vehicle platoon control problem, or so-called collaborative adaptive cruise control (CACC) \cite{varaiya1993smart}, has been widely studied in the literature and several solutions have been proposed \cite{iihoshi2000vehicle,hedrick1991longitudinal,kato2002vehicle} and is a natural extension of ACC that leverages vehicular ad hoc networks and vehicle to vehicle (V2V) communication. This problem has been well studied in the context of MPC control strategies \cite{schmied2015nonlinear,stanger2013model,sancar2014mpc}, which have a natural advantage of using the predictive nature of MPC and then sharing these predictions over the wireless channel, in order to improve overall system performance.

However, these works tend to make the same assumptions about wireless communication mentioned above. These assumptions generally fall into one of two categories: (1) perfect communication, i.e. at every time step of the model predictive controller, the vehicle gets new trajectory predictions from vehicles in the platoon; or (2) imperfect communication where the performance is stationary, and so a vehicle will receive a new packet with a known (stationary) probability.

Environmental conditions such as a bridge overpass, buildings, or other vehicles impact the quality of the communication channel. Many of these conditions involve multipath reflection, for example an adjacent semi-truck might cause multipath interference between two communicating vehicles.
In addition, the channel is affected by the {\em trajectory of the vehicle itself}; moving the transceivers closer or farther apart might positively (or negatively) impact multipath reflection.

This paper leverages recent and ongoing advances in wireless channel predictive capability, where knowledge of the scene can be used to predict what will happen to the channel \cite{kim2006accurate,de2005high,gnawali2009collection,farkas2006pattern}. The related literature uses ``black box'' or ``grey box'' approaches \cite{liu2011foresee,boano2010triangle}, where the function used to predict wireless properties is unknown (and possibly unknowable). State of the art techniques use some form of artificial neural network, or ANN \cite{elnaggarbayesian}. The notion of a black box is common in the field of machine learning, particularly those techniques that use ANNs, but this presents a significant challenge in the context of safety broadly, and for MPC specifically.

In the following work, we leverage recent advances in data-driven MPC \cite{rosolia2018learning,brunner2017repetitive,rosolia2017robust} that learn from previous iterations of a control task and provide guarantees on safety and improved performance at each iteration. 
In particular, we introduce a formulation of MPC for vehicle platooning that accounts for imperfect communication, and then design a LMPC control scheme that leverages the notion of predictive capability for wireless channel quality, in order to obtain better platoon performance.
The contributions of this paper are: 
\begin{enumerate}
    \item formulation of a (L)MPC problem that can handle decision variables or objective functions that derive from an unknown or unknowable function, for example variables that are generated by an artificial neural net, and
    \item extension of LMPC, adapted to handle dynamic environments and/or time-evolving constraints, in a computationally tractable manner.
\end{enumerate}

The paper is organized as follows: section \ref{sec:prelim} provides preliminaries about LMPC and vehicle dynamics,
and then section \ref{sec:problem} presents the MPC model for vehicle platooning in which the communication channel is considered to be perfect and the packet is delivered to the following vehicle without any delay. 
Section \ref{sec:lmpc-cavs} then provides a modified LMPC for CAVs in three parts: section \ref{sub:sr-lmpc} shows a tailored LMPC for dynamic environments, section \ref{sub:SR-LMPC with Uncertainty} takes into account the uncertainty of the communication channel, and section \ref{sub:MINLP-to-NP} describes how the presented Mixed Integer Nonlinear Problem can be transformed to a Nonlinear Problem. In section \ref{sec:example} the simulation results of SR-LMPC for two connected vehicles are presented, and section \ref{sec:conclusions} makes concluding remarks.

\section{Preliminaries}\label{sec:prelim}

This section is based on the original work of \cite{rosolia2018learning}. Beginning with a discrete time system
\[
x_{t+1} = f\left(x_t,u_t \right), \label{eq:nldynamics}
\]
where $x\in\mathbb{R}^n$ and $u\in\mathbb{R}^m$ are the system state and input, respectively, assume that $f(\cdot ,\cdot)$ is continuous and that state and inputs are subject to the constraints
\[x_t\in \mathcal{X}, \quad u_t \in \mathcal{U} \quad \forall t\ge0.\]\label{eq:constraints}
LMPC solves the following infinite horizon optimal control problem iteratively:
\begin{subequations}\label{eq:generic-objective}
    \begin{align}
        J_{0\to\infty}^* &= \min_{u_0,u_1,\dots}\sum_{k=0}^{\infty} h\left(x_k , u_k\right)\tag{\ref{eq:generic-objective}}\\
        \text{s.t.}\ & x_{t+1} = f\left(x_t,u_t \right) \quad \forall k \ge 0 \label{eq:dynamic-const}\\
        & x_0 = x_S \label{eq:init} \\
        & x_k \in \mathcal{X}, \quad u_k \in \mathcal{U}\quad \forall k \ge 0 \label{eq:mpc-xu-consts}
    \end{align}
\end{subequations}
where equations (\ref{eq:dynamic-const}) and (\ref{eq:init}) represent the system dynamics and the initial condition, and (\ref{eq:mpc-xu-consts}) are the state and input constraints. LMPC assumes that the stage cost $h(\cdot,\cdot)$ in equation (\ref{eq:generic-objective}) is continuous and satisfies
\begin{multline}
    h\left(x_F,0\right)=0 \ \text{and} \ h\left(x_t^j , u_t^j\right)\succ 0 \ \forall x_t^j \in \mathbb{R}^n \setminus \left\{x_F\right\},\\ u_t^j \in \mathbb{R}^m \setminus \left\{0\right\} \label{eq:final-state-cond}
\end{multline}
where the final state $x_F$ is a feasible equilibrium for the unforced system (\ref{eq:nldynamics})
\[f(x_F,0)=x_F.\]
At the $j^\text{th}$ iteration of LMPC, the vectors
\begin{subequations}\label{eq:jth-iteration}
    \begin{align}
      \textbf{u}^j &= \left[ u_0^j, u_1^j, \dots , u_t^j, \dots \right] \label{eq:input1}\\
      \textbf{x}^j &= \left[ x_0^j, x_1^j, \dots , x_t^j, \dots \right]\label{eq:state1} 
    \end{align}
\end{subequations}
collect the inputs applied to system (\ref{eq:nldynamics}) and the corresponding state evolution. In (\ref{eq:jth-iteration}), $x_t^j$ and $u_t^j$ denote the system state and the control input at time $t$ of the $j^\text{th}$ iteration. We assume that at each $j^\text{th}$ iteration, the closed loop trajectories start from the same initial state
\[x_0^j=x_{S},\quad \forall j \ge 0.\]\label{eq:}

\subsection{Sampled Safe Set}
A key notion of LMPC is that it exploits the iterative nature of control design. For every $k^{th}$ iteration that successfully steers the system to the terminal point $x_F$ (i.e.,$\forall k:\lim_{t\to\infty} x_t^k=x_F $), the trajectory $\textbf{x}^k$ is a subset of {\em sampled Safe Set} $\mathcal{SS}^j$, which is defined as:
\[
\mathcal{SS}^j = \left\{ \bigcup_{i\in M^j} \bigcup_{t=0}^{\infty} x_t^i\right\}
\]
where
\[
M^j = \left\{k\in [0,j]:\lim_{t\to\infty} x_t^k = x_F \right\}\label{eq:feasible-indices}
\]
$\mathcal{SS}^j$ is the collection of all state trajectories at iteration $i$ for $i\in M^j$. $M^j$ in (\ref{eq:feasible-indices}) is the set of indices $k$ associated with successful iterations of MPC for $k<j$. It follows that $\mathcal{SS}^i \subseteq \mathcal{SS}^j\ \forall i\le j$. $\mathcal{SS}^j$ is a subset of the maximal stabilizable set because, for every point in the set, there exists a feasible control action that satisfies the state constraints and steers the state toward $x_F$.

\subsection{Iteration Cost}
Define a function $Q^j(\cdot)$ over $\mathcal{SS}^j$ that assigns to every point in $\mathcal{SS}^j$ the minimum cost-to-go along the trajectories in sampled safe set,
\[
Q^j(x)=\begin{cases} \min\limits_{(i,t)\in F^j(x)} J_{t\to\infty}^i(x), & \text{if}\ x\in \mathcal{SS}^j \\
 +\infty & \text{if} \ x\notin  \mathcal{SS}^j\end{cases} 
\]
where $F^j(\cdot)$ is
\begin{multline}\label{eq:fj-assigner}
    F^j(x) = \Big\{ (i,t) : i\in [0,j], \ t\ge 0 \ \text{with } x_t^i=x;  \\
    \forall x_t^i\in \mathcal{SS}^j
\end{multline}
In other words, for every $x \in \mathcal{SS}^j$, $Q^j(x)$ not only assigns the optimal cost-to-go but also the pair $(i,t)$ that indicates the optimal iteration number in LMPC as well as the optimal time-to-go for that state.

\subsection{Properties of LMPC}
It can be shown \cite{rosolia2018learning} that, using the above notions of sampled safe set and iteration cost, the $j^\text{th}$ iteration cost is nonincreasing at each iteration and that the LMPC formulation is recursively feasible (state and input constraints at the next iteration are satisfied they are satisfied at the current iteration). It should also be noted that LMPC solves the infinite time optimal control problem by solving at time $t$ of iteration $j$ a finite time constrained optimal control problem. 

\section{Problem Formulation: MPC for CAVs}\label{sec:problem}
Consider that there are $\mathit{c+1}$ vehicles in the platoon, $i=\{0,1,\dots,c\}$, where $0$ refers to the lead car. The dynamical system of the $i^{th}$ vehicle is 
\begin{equation} \label{dynamical System}
    x_i(t+1)=f_i(x_i(t),u_i(t))
\end{equation}
where, $x_i(t)\in\mathbb{R}^n$ and $u_i(t)\in\mathbb{R}^m$ show the state and control input vectors of the $i^{\text{th}}$ vehicle at time step $t$, respectively. Also, assume that $f_i:\mathbb{R}^n\times\mathbb{R}^m\to\mathbb{R}^n$ is a smooth function that evolves the state of the vehicle $i$ through the time horizon, $N$.

The proposed control architecture is formulated as model predictive control that repeatedly solves the following finite horizon optimization problem:
\begin{subequations}\label{eq:mpc-cavs-generic}
    \begin{align}
         \min_{U_{t\to t+N|t}} & J = \sum_{k=t}^{t+N}
        \|{x}_i(k|t)-{x}_{i-1}(k|t)\|^2_{P^1_{i}} \notag\\
        & \quad +\|{x}_i(k|t)-{x}^{ref}(k)\|^2_{P^2_{i}} \tag{\ref{eq:mpc-cavs-generic}}\\ 
        & \quad + {\|u_i(k|t)\|^2_{P^3_{i}}} + Q_i(x_i(t+N|t)) \notag\\
        \text{s.t.}\ & x_{i}(k+1|t) = f_i\left(x_i(k|t),u_i(k|t) \right)  \label{eq:mpc-dynamic-const1}\\
        & x_i(t|t) = x_i(t) \label{eq:mpc-dynamic-const2} \\
        &\gamma_i(x_i(k|t),x_{i-1}(k|t))\ge \Gamma_{i}(k)  \label{eq:mpc-dynamic-const3}\\
        &\underline{x}_i \le x_i(k|t) \le \overline{x}_i \label{eq:mpc-dynamic-const4}\\
        &\underline{u}_i \le u_i(k|t) \le \overline{u}_i \label{eq:mpc-dynamic-const5}
         \\&\qquad \forall k \in \{t,\dots,t+N\} \notag
    \end{align}
\end{subequations}
where $x_i(k|t)$ and $u_i(k|t)$ are the state and control input at step $k$ predicted at time $t$, respectively. 


The objective function represents the tradeoff between tracking the preceding vehicle, reference tracking, and control effort with tuning matrices $P^1_i$, $P^2_i$, and $P^3_i$, respectively. $Q_i$ represents the terminal cost.
The first constraints (\ref{eq:mpc-dynamic-const1}) and (\ref{eq:mpc-dynamic-const2}) are the dynamical system of the vehicle and its initial state. The last two constraints (\ref{eq:mpc-dynamic-const4}) and (\ref{eq:mpc-dynamic-const5}) define the bounds on state variables and control inputs. The safety condition is enforced by inequality (\ref{eq:mpc-dynamic-const3}) in which $\Gamma_{i}(t)$ is \textit{Time To Collision} which is given as a parameter to the model: 
\begin{equation} \label{eq:ttc}
  \gamma_i(x_i(k),x_{i-1}(k)) =
  \begin{cases}
  \infty                 & v_{i-1}(k)\ge v_{i}(k)\\
  \frac{d_i(k)}{v_{i}(k)-v_{i-1}(k)} & v_{i-1}(k)<v_{i}(k)
  \end{cases}
\end{equation}
where $v_{i}(k)$ shows the velocity of vehicle $i$ and $d_i(k)$ is the euclidean distance between two vehicles $i-1$ and $i$ at time $k$ which can be calculated from their state vectors. 
The optimal solution to this problem is given by
\begin{equation}
  \label{eq:optimial_sol_MPC}
  \begin{gathered}
    \mathrm{\textbf{X}}_i(t)=\left[x_i(t+1|t), x_i(t+2|t),...,x_i(t+N|t)\right]\\
    \mathrm{\textbf{U}}_i(t)=[u_i(t|t), u_i(t+1|t),...,u_i(t+N-1|t)]
  \end{gathered}
\end{equation}
and a receding horizon control law applies the first control input $u_i(t|t)$  that transfers the system to the new state $x_i(t+1|t)$, and the process is repeated from $t+1$.

To solve this problem, it is necessary to assume that there is a feasible trajectory from the starting point $x_{i,s}$ to the goal point $x_{i,g}$, which includes the feasibility of the model at time $t$ (i.e. the solution space of the MPC is not empty).

Like many existing works \cite{firoozi2018safe}, the above MPC platooning control scheme assumes that the communication channel works perfectly, and that each vehicle has access to the (current) trajectory prediction of the preceding vehicle. We seek to relax this assumption. To take into account the uncertainty in the wireless channel, we adapt and extend LMPC for connected vehicles and platooning, which has previously been presented for a single vehicle in a static environment \cite{rosolia2018learning}.
In the following sections, for the purpose of simplicity, only one iteration of model (\ref{eq:mpc-cavs-generic}) is considered for developing the algorithm.


\section{Learning Model Predictive Control for CAVs}\label{sec:lmpc-cavs}
The original LMPC formulation presented above is designed for an entire trajectory from an initial state to a final state, in a static environment. As those authors acknowledge, LMPC is computationally expensive; in addition, LMPC cannot be applied in a dynamic environment in its original form. To overcome these limitations, we propose two concepts. The first notion involves encoding the dynamics of obstacles both in the constraints and the objective function of the optimization problem. In addition, LMPC is reformulated from an end-to-end planning problem (i.e. find an optimal trajectory from $x_S\to x_F$), which typically has better performance in terms of the number of required iterations and degree of optimality, as the planning horizon, $N$, increases \cite{rosolia2017learning}. In a dynamic context, a shorter planning horizon trades convergence to global optima with the ability to overcome computational issues with LMPC in general. 

\begin{figure}[h]
    \centering
    \hspace{-10pt}\includegraphics[width=1.04\linewidth]{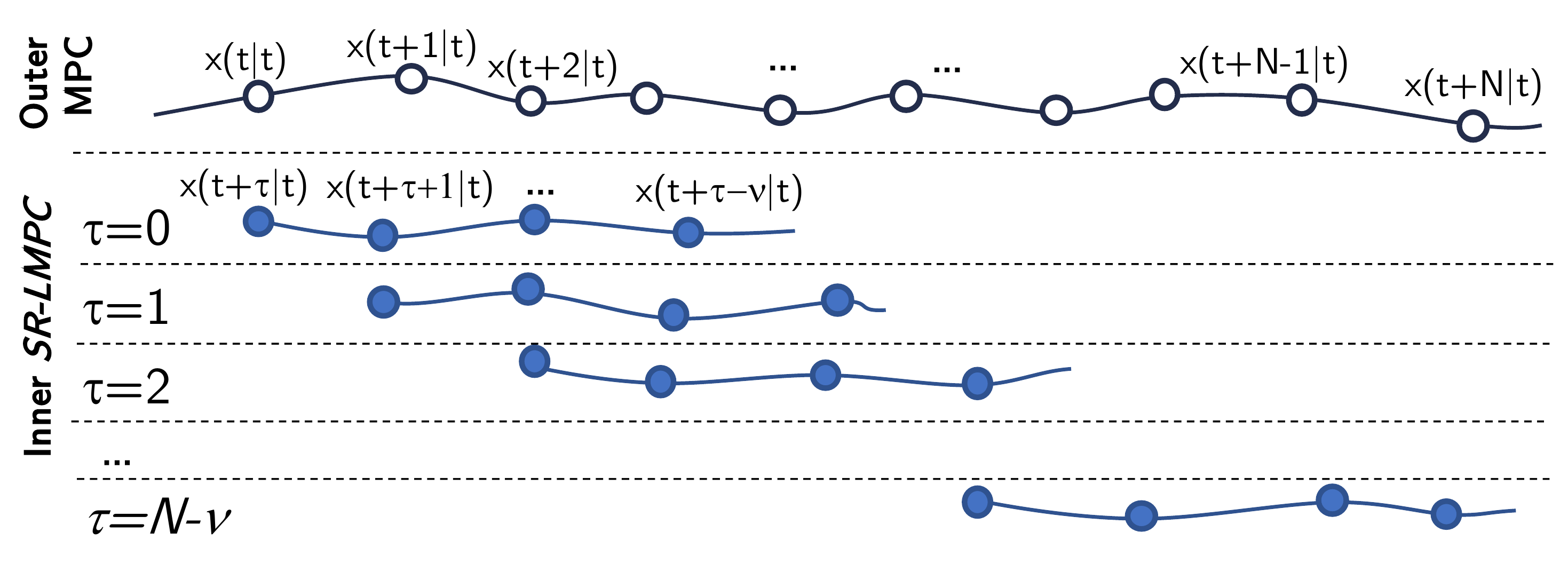}
    \caption{Workflow of SR-LMPC Control Architecture: for each $t$ in the nominal MPC outer loop, SR-LMPC computes $N-\nu$ shorter, receding-horizon runs over data-driven decision variables}
    \label{fig:SR-LMPC-scheme}
\end{figure}

In addition, LMPC is extended to account for decision variables with unknown dynamics or data-driven approximations. The control architecture includes a (nominal) outer loop, or high level, motion planner that generates a candidate set of feasible trajectories. LMPC works on the inner loop to converge to a dynamically feasible trajectory with optimal (or improved) performance over data-driven decision variables. The approach is depicted graphically in Fig. \ref{fig:SR-LMPC-scheme} and pseudocode for this architecture is shown at the end of this section in Algorithm \ref{alg2}.

\subsection{Short Range LMPC}\label{sub:sr-lmpc}
\begin{assumption}
Assume that the preceding vehicle, $i-1$, has converged on its own optimal trajectory and that this is known for the time horizon, $N$:  
\begin{equation}
  \label{eq:optimal_trajectory}
  \begin{gathered}
    \textbf{X}^{*}_{i-1}(t)=[x_{i-1}^{*}(t|t),x_{i-1}^{*}(t+1|t),\dots,x_{i-1}^{*}(t+N|t)]\\
  \end{gathered}
\end{equation}
\end{assumption}
This assumption should be relaxed in future work, but (under) approximates an autonomous system's ability to predict its own trajectory into the future and communicate this to the platoon. In section \ref{sub:SR-LMPC with Uncertainty}, the assumption on perfect communication will be relaxed and a different control strategy will be adopted.

\begin{assumption}
 SR-LMPC assumes the existence of a feasible trajectory from the current state of the ego vehicle (i.e. the follower), $x(t) \to x(t+N)$ at the first iteration but with no assumptions on optimality, given as
\begin{subequations}
    \begin{align}
    \textbf{X}^0(t)=[x^0(t+1|t), x^0(t+2|t), \dots, x^0(t+N|t)]\label{eq:init-state-dynlmpc}\\
    \textbf{U}^0(t)=[u^0(t|t), u^0(t+1|t),\dots,u^0(t+N-1|t)]\label{eq:init-input-dynlmpc}
    \end{align}
\end{subequations}
where $x^0(k|t)$ and $u^0(k|t)$ are the state and control input vectors of the ego system at time $k$ that have been calculated at time step $t$. The superscript shows the iteration number of LMPC, which starts from $0$ and is denoted by index $\ell=\{0,1,...,L\}$. SR-LMPC keeps only a record of successful iterations in this set, with the total number of successful iterations completed by the algorithm given by $L$. The information of each trajectory is saved in set $\mathcal{DS}$ if it is completed successfully, which implies no collisions (see (\ref{eq:ttc})) but not necessarily optimality. Given an arbitrary, feasible initial trajectory stored in $\mathcal{DS}^0$, the dynamic safe set is iteratively built as
\begin{multline}
    \mathcal{DS}^L=\Bigg\{\bigcup_{k=1}^N \Big\{x^L(t+k|t) \ \Big| \ x_k^L \in \mathcal{X}^{\mathrm{d}_t},
    \\ u_k^L \in \mathcal{U} \ \forall k \in \{t,t+1,\dots,t+N\} \Big\}\Bigg\}\cup  \mathcal{DS}^{L-1} \label{eq:dynamic-safe-set}
\end{multline}
where $x\in\mathcal{X}^{\mathrm{d}_t}$ represents the dynamic constraints on state imposed by the time-to-collision for the predicted leader trajectory at time $t$. 
\end{assumption}

The cost-to-go of state $x(k|t)$ is denoted by $q(x(k|t)|x_{i-1}^*(k\to N|t))$ and is defined as ``the trajectory cost of the ego system from $x(k|t)$ to $x(k+N|t)$ given that the states of the leading vehicle are $x_{i-1}^*(k|t),...,x_{i-1}^*(k+N|t)$". A backward calculation is applied to find the cost-to-go for each state in set $\mathcal{DS}^{L}$. To start, $q_i^{\ell}(x_i^{\ell}(N|t)|x_{i-1}^*(N|t))$ can be approximated by path planning algorithms \cite{jafarzadeh2018exact} or simply assumed 1-norm distance between two vehicles at time $N$.
\begin{multline}
  \label{eq:cost_to_go}
  q_i^{\ell}(x_i^{\ell}(k|t)|x_{i-1}^*(k \to N|t))=z_i(x_i^{\ell}(k|t),u_i^{\ell}(k|t))\\ +q_i^{\ell}(x_i^{\ell}(k+1|t)|x_{i-1}^*(k+1 \to N|t))
\end{multline}
for $k=\{N-1, ... , t\}$, and where $z_i(x_i^{\ell}(k|t),u_i^{\ell}(k|t))$ is the stage cost and is defined as the cost of control effort to transfer the state of the system from $x_i^{\ell}(k|t)$ to $x_i^{\ell}(k+1|t)$ by applying input $u_i^{\ell}(k|t)$ at iteration $\ell$ and time $k$. The overall performance of the controller at iteration $\ell$ occurs when $k=t$.

The SR-LMPC formulation for vehicles in a platoon takes the form of the following constrained mixed integer optimization problem (MINLP). SR-LMPC uses the generic form of problem (\ref{eq:mpc-cavs-generic}), with two notable exceptions. 
First the time horizon is modified, and significantly for the purposes of computation shortened: $\{t+\tau,t+\tau+1,...,t+\tau+\nu\}$, where $\nu <N$ is the SR-LMPC time horizon and $t+\tau$ is the starting time, for all $\tau =\{0,1,..., N-\nu\}$.

The objective function is modified to
\begin{multline}
    J_{t+\tau \to t+\tau+\nu}^L(x_i^L(t+\tau))= 
    \\ \sum_{k=t +\tau}^{t+\tau+\nu-1} \|x_i(k+1|t)-x_{i-1}^*(k+1|t)\|^2_{P^1_{i}}+ 
     \|u_i(k|t)\|^2_{P^2_{i}}
    \\ +\sum_{\ell=0}^{L-1}{\sum_{\eta=0}^N{\mathrm{\zeta}_i^{\ell}(\eta) q_i^{\ell}(x_i^{\ell}(t+\eta|t)|x_{i-1}^*(\eta \to N|t))}}     \label{LMPC:obj}
\end{multline}
All of the constraints (\ref{eq:mpc-dynamic-const1})-(\ref{eq:mpc-dynamic-const5}) hold, given appropriate sets $k\in\{t+\tau,\dots,t+\tau-\nu\}$. 

Finally, SR-LMPC then adds the following constraints
\begin{subequations}\label{eq:SR-LMPC-newconstraints}
\begin{align}
    &x_i(t+\tau+\nu|t)=\sum_{\ell=0}^{L-1} {\sum_{\eta=0}^N {\zeta_i^{\ell}(\eta) x_i^{\ell}(t+\eta|t)}} \label{LMPC:con4}\\
    &\sum_{\ell=0}^{L-1}{\sum_{\eta=0}^N {\zeta_i^{\ell}(\eta)}}=1 \label{LMPC:con5}\\
    & \mathrm{\zeta}_i^{\ell}(\eta)\in\{0,1\}, \forall\ell=\{0,..,L-1\}, \forall \eta=\{0,...,N\} \label{LMPC:con6}
    \end{align}
\end{subequations}

The following vectors show the optimal solution for this model:
\begin{equation}
  \label{eq:SR-LMPC}
  \begin{gathered}
    \mathpzc{X}_i(t+\tau)=[\mathpzc{x}_i(t+\tau+1),\mathpzc{x}_i(t+\tau+2),...,\mathpzc{x}_i(t+\tau +\nu)]\\
    \mathpzc{U}_i(t+\tau)=[\mathpzc{u}_i(t+\tau),\mathpzc{u}_i(t+\tau +1),...,\mathpzc{u}_i(t+\tau +\nu-1)]
  \end{gathered}
\end{equation}

This model assigns a binary decision variable, $\zeta_i^{\ell}(\eta)$, to each of the states in $\mathcal{DS}^L$ (\ref{eq:dynamic-safe-set}) and selects one of them as the terminal state, (\ref{LMPC:con4}). Because only one of these states can be chosen as the terminal state $x_i(t+N|t)$, the summation of the binary variables should be exactly one, as shown in (\ref{LMPC:con5}). The last term of the objective function (\ref{LMPC:obj}) determines the best value for cost-to-go $q_i^{\ell}(x_i^{\ell}(t+\eta|t)|x_{i-1}^*(\eta \to N|t))$ and its associated state. Assuming that the optimal state in set $\mathcal{DS}^L$ is $x_i^{\ell^{\ast}}(t+\eta^{\ast}|x_{i-1}^*(\eta^{\ast} \to N|t))$, the solution for system state is the vector
\begin{equation}
  \label{eq:SRLMPC_best_terminal_state}
  \begin{gathered}
    \mathpzc{x}_i(t+\tau+\nu)=x_i^{\ell^{\ast}}(t+\eta^{\ast}|x_{i-1}^*(\eta^{\ast} \to N|t))\\
  \end{gathered}
\end{equation}

After finding the optimal solution for SR-LMPC as an MINLP model, the first step of the control input vector, $\mathpzc{U}_i(t+\tau)$, is implemented and the related state and control vectors are saved in the trajectory of iteration $L$:
\begin{equation}
  \label{eq:aplly_the_first_step}
  \begin{gathered}
    x_i^L(t+\tau +1|t)=\mathpzc{x}_i(t+\tau+1)\\
    u_i^L(t+\tau|t)=\mathpzc{u}_i(t+\tau)
  \end{gathered}
\end{equation}

The updated trajectory in the current iteration, $L$, is as follows (notice the slightly different symbol for $\mathrm{x}$ and $\mathrm{u}$, indicating the first step of the receding horizon LMPC along the time-shift $\tau$):
\begin{subequations}
  \label{eq:LMPC_trajectory}
  \begin{align}
    \mathrm{\textbf{X}}_i^L(t)&=[x_i^L(t+1|t), x_i^L(t+2|t),\dots,x_i^L(t+\tau+1|t)]\notag\\
    \mathrm{\textbf{U}}_i^L(t)&=[u_i^L(t|t),u_i^L(t+1|t),\dots,u_i^L(t+\tau|t)]\tag{\ref{eq:LMPC_trajectory}}
  \end{align}
\end{subequations}
for $\tau =\{0,1,...,N-\nu\}$.

\subsection{SR-LMPC with Uncertain Communication Channel} \label{sub:SR-LMPC with Uncertainty}
\begin{figure*}[tp]
    \centering
    \hspace{-10pt}\includegraphics[width=1\linewidth]{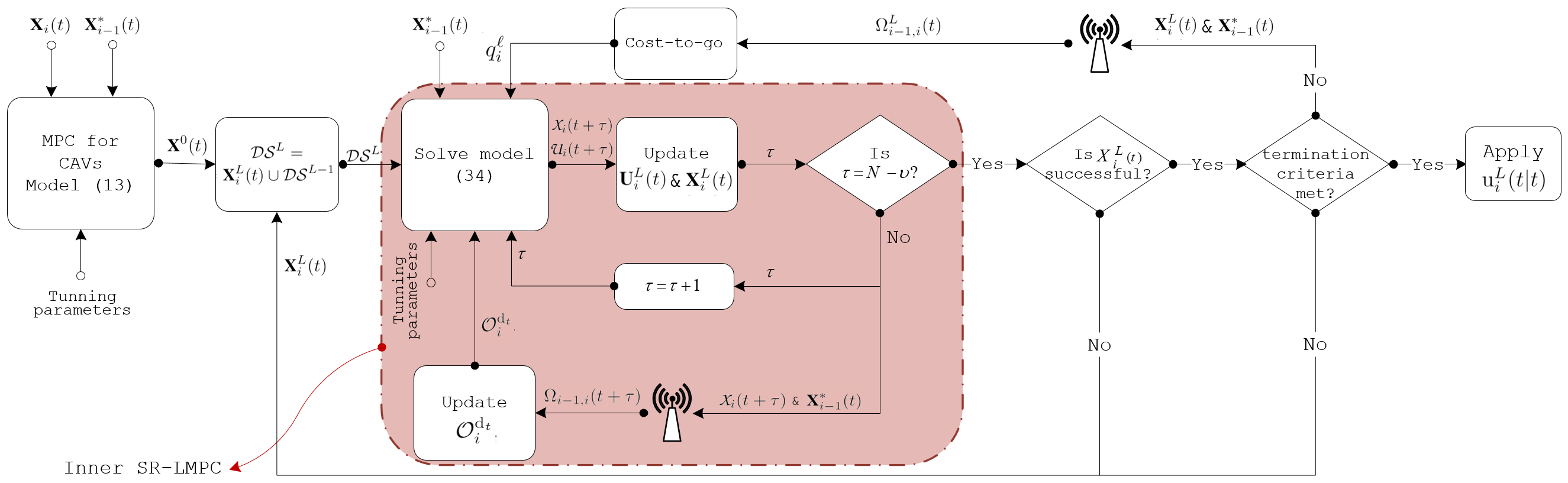}
    \caption{Block diagram of SR-LMPC}
    \label{fig:SR-LMPC-blockDiagram}
\end{figure*}

The formulation in \ref{sub:sr-lmpc} assumes no delay in the communication channel and the following vehicle receives the lead vehicle's motion plan for the current time horizon. 
However, the communication channel can be affected by numerous factors such as relative states of the communicating vehicles, adjacent buildings and vehicles, and so on. The quality of the communication channel is defined by Packet Delivery Rate, $PDR_t$, and its inverse at time $t$ provides an estimate of delivery time for each packet $\omega_t^{\ell}=\frac{1}{PDR_t^{\ell}}$. This estimate can be provided by an unknown function, e.g. an artificial neural network like an LSTM~\cite{elnaggarbayesian}. The radio icons in Fig.~\ref{fig:SR-LMPC-blockDiagram} represent the LSTM's mapping from vehicle state to channel estimation.

At iteration $\ell$ and time step $t$, the delivery time in range $ j\cdot dt \leqslant \omega_t^{\ell} < (j+1)dt$, where $j \in \mathbb{Z}^{\geqslant}$, and $dt$ is the time constant used for the discretized dynamics in the general MPC formulation. 
This model of communication implies that the vehicle $i$ should rely on the most recent packet available at time $t-dt$ 
to calculate its motion policy. Assume that vehicle $i$ has received several packets at time step $t$ from the $i-1$ lead vehicle, which contain the trajectory predictions from whatever time stamp they were sent:   $\{\mathrm{\textbf{X}}_{i-1}^*(t_1),\mathrm{\textbf{X}}_{i-1}^*(t_2),,...,\mathrm{\textbf{X}}_{i-1}^*(t_m)\}$, where $t_1 < t_2 < ... < t_m \leqslant t $. Because the packet containing $\mathrm{\textbf{X}}_{i-1}^*(t_m)$ is the most up-to-date (though not necessarily {\em current}), it is used to calculate the motion policy of the following vehicle. 
\begin{multline}
  \label{eq:most_recent_packet}
    \mathrm{\textbf{X}}_{i-1}^*(t_m)=[x_{i-1}^*(t_m|t_m),x_{i-1}^*(t_m+1|t_m),\\
    \dots,x_{i-1}^*(t_m+N|t_m)]
\end{multline}

However, the first $t-t_m$  states are stale, i.e. the leader has already executed these states, and should be deleted from the packet. Therefore, the packet shrinks to:
\begin{multline}
  \label{eq:pruned_packet}
    \mathrm{\textbf{X}}_{i-1}^*(t_m)=[x_{i-1}^*(t|t_m),x_{i-1}^*(t+1|t_m),\\ \dots,x_{i-1}^*(t_m+N|t_m)]
\end{multline}
The length of this packet determines the real time horizon that the following vehicle can consider:
\begin{equation}
  \label{eq:changing time horizon}
  \begin{gathered}
  \mathrm{N}_i(t)=N-(t-t_m)\\
  \end{gathered}
\end{equation}
where, $\mathrm{N}_i(t)$ denotes the length of time horizon at time $t$, $\mathrm{N}_i(t) \leqslant N $. Fewer useful states may cause vehicle $i$ to take more conservative policies to satisfy the constraints, and this results in more cost for the entire whole trajectory; i.e. the controller might find local optima due to lack of longer-term information (see results in section \ref{sec:example}). 

One of the key innovations of our approach is adding awareness of the cost to physical system performance due to communication delays. Therefore, the updated formulation includes communication delay, $\omega_t^{\ell}$, in the objective function to penalize the trajectory that results in low quality (predicted) communication channel. Thus, the optimality condition forces the model to generate a trajectory in which the quality of the communication achieves an acceptable level, to avoid extra cost in the objective while satisfying the constraints. The vector of packet delivery time for whole time horizon is 
\begin{multline}
      \label{eq:packet_delivery_time_MPC}
    {\Omega}_{i-1,i}^L(t)=[{\omega}_{i-1,i}^L(t+1|t), {\omega}_{i-1,i}^L(t+2|t),\\
    \dots,{\omega}_{i-1,i}^L(t+N|t)]
\end{multline}

We assume that this vector is calculated and given by another system that is outside the scope of this work.
One major challenge is that the delivery time of the packets cannot be calculated before determining the trajectory, as communication depends on the relative location of vehicles and characteristics of their surrounding environment in different time steps. Also, after generating a trajectory, $\mathrm{\Omega}_{i-1,i}^L(t)$ is calculated outside of the control algorithm (e.g. through methods referenced in section \ref{sec:intro}), and because the model has not considered this variable in finding the solution, $\mathrm{N}_i(t)$ does not necessarily have an optimal value in following time steps. This can result in increases in the cost function.

To solve this problem, we exploit the repetitive nature of LMPC to consider $\mathrm{\Omega}_{i-1,i}^L(t)$ as a data-driven decision variable in the optimization process. This is done by two modifications in the model. First, we update $q_i^{\ell}$ in the objective function to include the cost associated with the delay in the communication channel if a specific terminal state is chosen. Observe that the cost-to-go vector is updated after a complete trajectory is generated. Second, a dynamic constraint is added to the model to represent the area where the communication loss occurs and is updated at each inner SR-LMPC iteration. For the first part, we define the following cost function
\begin{equation}
\label{eq:cost_function_packet_delivery_time_SR-LMPC}
    \sum_{k=t}^{t+N} \alpha^{(t+N)-k}{\omega}_{i-1,i}(k), 0<\alpha<1
\end{equation}
If $\alpha$ is close to zero, the cost function will assign more weight to the most recent values, but if it is close to one, the weight of delivery time at all time steps will receive almost equal weight. Now the cost-to-go in the original objective function (\ref{LMPC:obj}) can be reformulated as:
\begin{align}
  q_i^{\ell}(x_i^{\ell}&(k|t)| x_{i-1}^*(k \to N|t)) = \notag\\ & z_i(x_i^{\ell}(k|t),u_i^{\ell}(k|t)) + 
   {\omega}_{i-1,i}(x_i^{\ell}(k|t)|x_{i-1}^*(k|t)) \notag\\
   &+ (\alpha -1)\sum_{j=k+1}^{t+N}\alpha^{j-k}{\omega}_{i-1,i}(j)  \notag\\
   &+ q_i^{\ell}(x_i^{\ell}(k+1|t)|x_{i-1}^*(k+1 \to N|t))\notag\\
   =& \sum_{j=k}^{t+N}\left[ \alpha^{j-k}{\omega}_{i-1,i}(j) + z_i(x_i^{\ell}(j|t),u_i^{\ell}(j|t)) \right]
\end{align}

The updated cost-to-go contains the cost of delivery time of communication channel and by replacing it in the objective function, the model will chose a terminal state that has lower cost in terms of communication channel, $\mathrm{\Omega}_{i-1,i}^L(t)$, and stage costs, $z_i$.For the second part, a boundary on the state vector is defined as a dynamic constraint in the model
\begin{equation}
  \label{eq:dead-zone obstacle}
  \begin{gathered}
  x_i(k|t) \notin \mathcal{O}^{\mathrm{d}_t}_i \\
  \end{gathered}
\end{equation}

Note that this is a time-variant constraint and it depends on the states of the lead and ego vehicles. $\mathcal{O}^{\mathrm{d}_t}_i$ is updated based on the vector of packet delivery time for SR-LMPC which is shown by
\begin{multline}
  \label{eq:packet_delivery_time_SR-LMPC}
    {\Omega}_{i-1,i}(t+\tau)=[{\omega}_{i-1,i}(t+\tau+1), {\omega}_{i-1,i}(t+\tau+2),\\
    \dots, {\omega}_{i-1,i}(t+\tau+\nu)]\\
\end{multline}

For simplicity, assume that at each iteration of Branch and Bound relaxation, the algorithm solves a convex quadratic model. Using the Interior Point Method, the computational complexity of finding $\epsilon-$scale optimum for a quadratic model is polynomial in the size of model ($n'$) and required accuracy ($\epsilon$), i.e. $O(n'log1/\epsilon)$ \cite{ye1989extension} .However, the worst-case number of iterations of B\&B algorithm is exponential $O(2^{(L-1)N})$, where $(L-1)N$ is the number of binary variables assigned to the vector $Q_i$. The size of model with time horizon $N$ is $(n+M)N$, resulting in computational complexity of $O(2^{(L-1)N}(n+m)Nlog1/\epsilon)$. The exponential part is dominant and yields in $O(2^{LN})$. On the other hand, the proposed method with smaller time horizon $\nu$ has time complexity of $O(2^{(L-1)N}(n+m)\nu(N-\nu)log1/\epsilon)$ at each outer loop iteration, which again yields $O(2^{LN})$. Then, the time complexity of the algorithm increases exponentially in the number of outer loop iterations, $L$. Without impacting computational complexity, shortening the time horizon from $N$ to $\nu$ enables the algorithm to utilize data at smaller but more frequent steps (inner SR-LMPC iterations) to improve the current trajectory by exploring the solution space with greater coverage. This approach noticeably decreases the number of overall trajectories, $L$, needed to converge to the optimum.

\subsection{Converting MINLP to Nonlinear Problem} \label{sub:MINLP-to-NP}

The formulation presented in \ref{sub:SR-LMPC with Uncertainty} computationally expensive to solve. 
The binary variables $\mathrm{\zeta}_i^{\ell}(\eta)$ and nonlinear constraints make the model Mixed Integer Nonlinear Programming (MINLP), which should be avoided in pursuit of algorithms that scale and/or can be applied in real time. Translating the above MINLP problem into the following Nonlinear Programming (NLP) formulation makes computing the optimal control policy more tractable. The control algorithm for SR-LMPC solves the following problem
\begin{subequations}\label{eq:NLP-problem}
    \begin{align}
        \min_{\mathpzc{U}_i(t+\tau)} & J_{t+\tau \to t+\tau+\nu}^L(x_i^L(t+\tau))= \tag{\ref{eq:NLP-problem}} \\ 
        & \sum_{k=t +\tau + 1}^{t+\tau+\nu}  \|{x}_i(k|t)-{x}_{i-1}^*(k|t)\|^2_{P^1_{i}} + \|u_i(k-1|t)\|^2_{P^2_{i}}\notag\\
        & +\sum_{\ell=0}^{L-1}{\sum_{\eta=0}^{\mathrm{N}_i(t)}{\mathrm{\zeta}_i^{\ell}(\eta) q_i^{\ell}(x_i^{\ell}(t+\eta|t)|x_{i-1}^*(\eta \to {\mathrm{N}_i(t)}|t))}}   \notag\\
        \text{s.t. } & x_i(k+1)=f_i(x_i(k),u_i(k)) \label{SR-LMPC:con1} \\
        & x_i(t+\tau|t)=x_i^L(t+\tau) \label{SR-LMPC:con2} \\
        & x_i(k|t) \in \mathcal{X}_i^{\mathrm{d}_t} \setminus \mathcal{O}^{\mathrm{d}_t}_i, u_i(k|t) \in \mathcal{U}_i \\ 
        & x_i(t+\tau+\nu|t)=\sum_{\ell=0}^{L-1} {\sum_{\eta=0}^{\mathrm{N}_i(t)} {\mathrm{\zeta}_i^{\ell}(\eta) x_i^{\ell}(t+\eta|t)}} \label{SR-LMPC:con4}\\
        & \sum_{\ell=0}^{L-1} {\sum_{\eta=0}^{\mathrm{N}_i(t)} {\mathrm{\zeta}_i^{\ell}(\eta)}}=1 \label{SR-LMPC:con5} \\
        & \mathrm{\zeta}_i^{\ell}(\eta)(1-\mathrm{\zeta}_i^{\ell}(\eta'))=0, \forall \eta' \geqslant \eta \label{SR-LMPC:con6_new} \\
        & 0 \leqslant \mathrm{\zeta}_i^{\ell}(\eta) \leqslant 1, \forall\ell=\{0,..,L-1\}, \notag \\ 
        & \qquad \forall \eta=\{0,...,{\mathrm{N}_i(t)}\} \label{SR-LMPC:con6}
    \end{align}
\end{subequations}

The difference in this model to the original SR-LMPC formulation (see constraints (\ref{LMPC:con4}) - (\ref{LMPC:con6})) is the values that $\mathrm{\zeta}_i^{\ell}(\eta)$ can take as a decision variable.Although constraint (\ref{SR-LMPC:con6}) allows $\mathrm{\zeta}_i^{\ell}(\eta)$ to take values from zero to one, constraint (\ref{SR-LMPC:con6_new}) limits them to be just one or zero. On the other hand, constraint (\ref{SR-LMPC:con5}) enforces them to be all zero except for one that should be valued one.

The overall SR-LMPC approach is shown in Algorithm \ref{alg2}. 

\begin{algorithm}
\caption{ SR-LMPC}
\label{alg2}
\begin{algorithmic}[1]
\STATE {{Solve model (\ref{eq:mpc-cavs-generic})}}
\STATE {Replace ${\textbf{x}}_i(t)$ in $\mathcal{DS}^0$}
\STATE {Set $L =0$ and $convergence = 0$}
\WHILE{$convergence = 0$}
\STATE {Calculate the prediction for ${\Omega}_{i-1,i}^L(t)$}
\STATE {Calculate the cost-to-go vector $q_i^L$}
\STATE {Set $L =L+1$}
\FOR{$\tau=1:N-\nu$}
\STATE {Update $ \mathcal{O}^{\mathrm{d}_t}_i$}
\STATE {Solve model (\ref{eq:NLP-problem})}
\STATE {Update trajectory $\mathrm{\textbf{X}}_i^L(t)$ and $\mathrm{\textbf{U}}_i^L(t)$}
\STATE {Apply $\mathpzc{u}_i(t+\tau)$}
\ENDFOR
\STATE {Update $\mathcal{DS}^{L}=\mathrm{\textbf{X}}_i^L(t)\cup  \mathcal{DS}^{L-1}$}
\IF{termination criteria are met}
\STATE {Set $convergence = 1$}
\ENDIF
\ENDWHILE
\STATE {Apply $\mathrm{u}_i^L(t|t)$}
\end{algorithmic}
\end{algorithm}


\section{Example Scenario - Uncertain Communication}\label{sec:example}
The application of a controller that adapts the vehicle’s longitudinal velocity based on the other vehicle’s states, based on uncertain communication over a V2V network, is not restricted to the car-following or platooning scenario. Similar formulations of the platooning or leader-follower control scheme from section \ref{sec:lmpc-cavs} can be extended to other applications like autonomous intersection control and improve either decentralized or centralized approaches~\cite{bashiri2017platoon,bashiri2018paim}.

However, the scenario developed for this paper focuses on a leader-follower scenario, where communication performance is influenced by a bridge overpass, which is known to have a negative effect \cite{schneider2016directional,cheng2013roadside,Hs2006}. The leading vehicle follows a simple trajectory with constant velocity. The following vehicle has stable initial conditions, i.e. before entering the bridge area it has converged to an optimal trajectory in terms of the tradeoff between following distance, control cost, and communication (we assume that communication performance is stable outside of the bridge scenario).

As the vehicles approach the bridge, the channel will begin to experience packet drops if the relative states achieve certain characteristics, for example combinations of relative speed and distance, due to multi-path interference with the bridge surface. It is also assumed that this deterioration in channel performance can be accurately predicted over time horizon $N$. In this simplified scenario, it is assumed that communication will drop out if both vehicles are in the so-called ``dead zone'' at the same time. One of the terms formulated in the objective function of the ego vehicle is avoiding states that cause communication loss (not entering into the dead zone) while satisfying other terms. Note that even in this simple scenario, SR-LMPC does not have access to this information in terms of a closed-form function in either the objective or constraints. Rather, at each (predicted or current) time step, an offline data-driven function is queried and the predictions are returned. 

Finally, the vehicle dynamics are formulated as a point mass system such that $\Dot{x}(t)=u(t)$, subject to input saturation. Each of these assumptions should be relaxed in future work, such as imperfect prediction accuracy and Dubins path dynamics. Relevant parameters of the scenario and model are shown in Table~\ref{tab:scenario_parameters}.
\vskip-13pt
\begin{table}[h]
    \centering
    \caption{Scenario and Model Parameters}\label{tab:scenario_parameters}
    \begin{tabular}{cccc} \hline
        $\underline{u}$,$\overline{u}$ & control input limits & m/s$^2$ & [-6,6] \\
        $N$ & time horizon of outer loop MP & -- & 70\\
        $\nu$ & time horizon of SR-LMPC & -- & 60\\
        $dt$ & sample time or discretization & s & 0.2\\
        $v_{\ell}^{ss}$ & lead vehicle speed (constant) & m/s & 30\\
        $v_{f,0}$ & follower vehicle initial speed  & m/s & 35\\
        $\Omega_{d}$ & dropout zone of bridge & m & 435-480 \\ \hline
    \end{tabular}
\end{table}

\begin{figure}[h]
    \vskip-6pt
    \centering
    \includegraphics[width=1.03\linewidth,trim={0.7in 0.4in 0.8in 0.6in},clip]{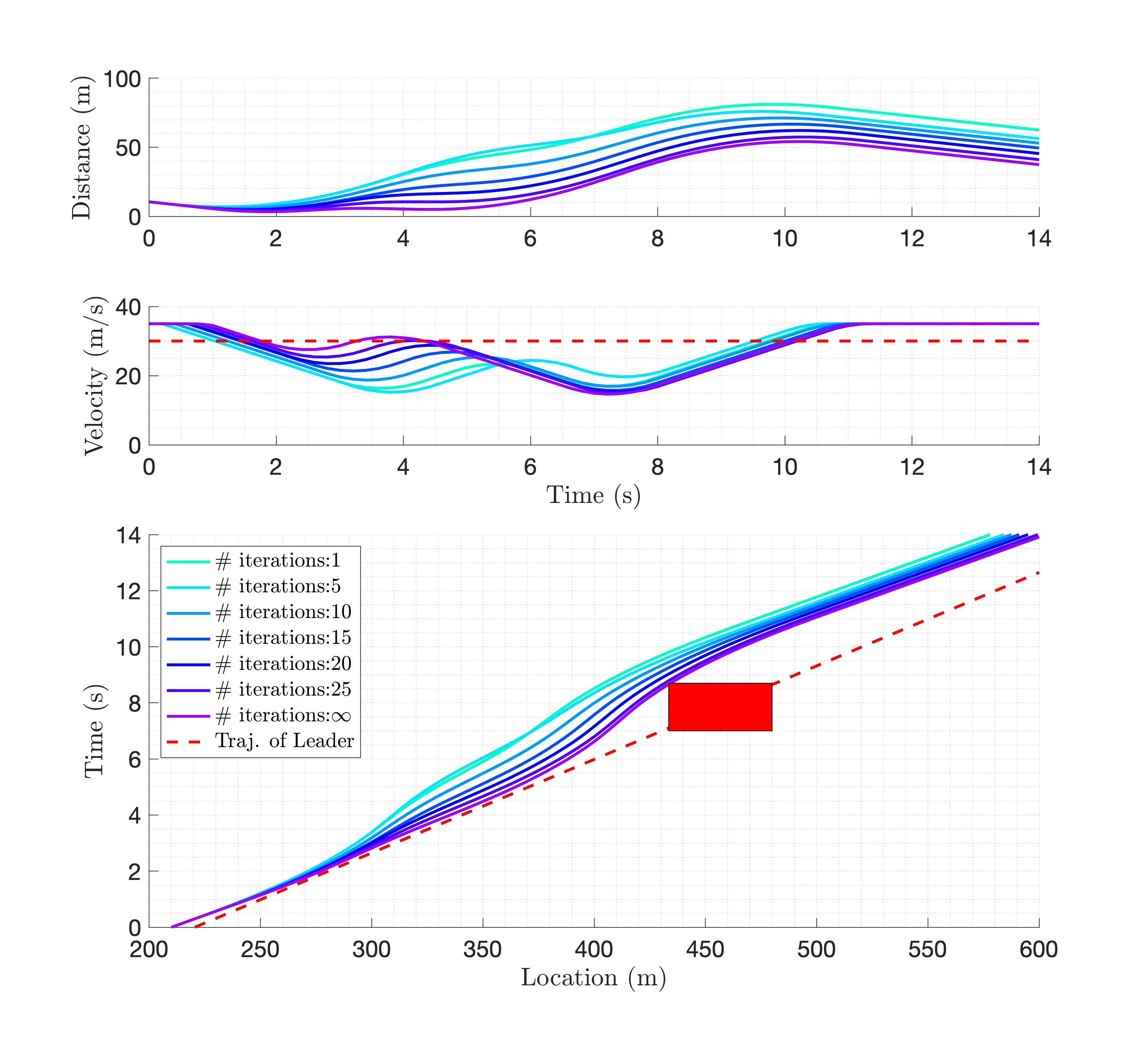}
    \vskip-11pt
    \caption{Car-following scenario: lead vehicle approaches bridge overpass where dropouts occur if both vehicles are in $[430,480]\ m$. LMPC performance with blackbox prediction of communication.} 
    \label{fig:lmpc-results-bridge}
\end{figure}

\begin{figure}[h]
    \vskip-5pt
    \centering
    \includegraphics[width=1\linewidth,trim={0.25in 0.1in 0.6in 0.3in},clip]{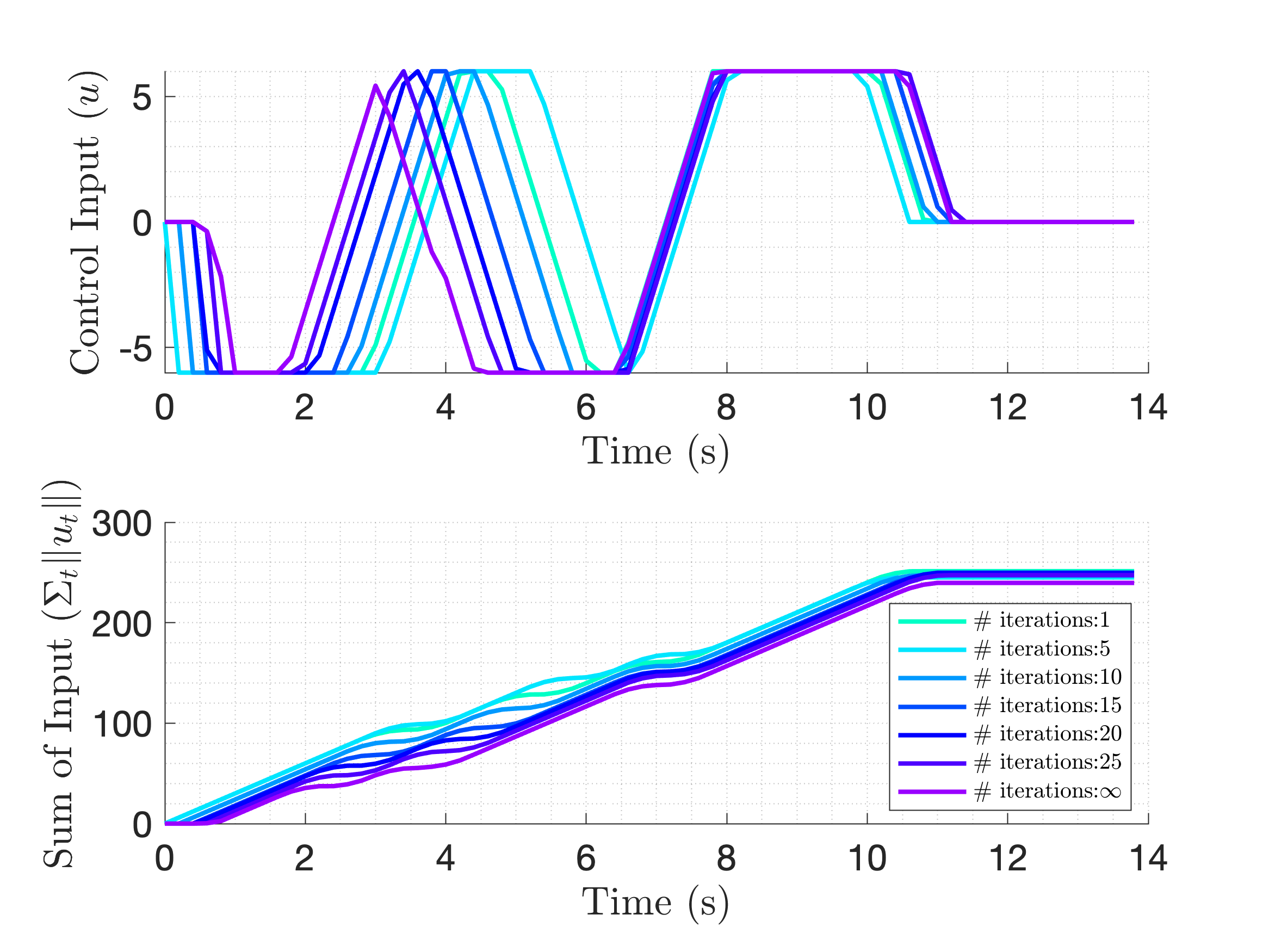}
    \vskip-10pt
    \setlength{\textfloatsep}{5pt}
    \caption{Control input for the car-following scenario. SR-LMPC converges to a solution that avoids input saturation and saves overall control cost.}
    \label{fig:control_input}
\end{figure}

\vskip-9pt
\noindent The optimization problems in (\ref{LMPC:obj}) and  (\ref{eq:NLP-problem}) are solved using CPLEX~\cite{cplex2009v12}. We seek to characterize whether, and to what extent, system performance can be improved with a ``look ahead'' function that can predict channel performance by means of a black box method. Results are shown in Fig.~\ref{fig:lmpc-results-bridge} and \ref{fig:control_input}.


It should be noted that the nominal scheme without SR-LMPC under communication uncertainty and/or prediction (section \ref{sub:sr-lmpc} and \ref{sub:SR-LMPC with Uncertainty}) proceeds agressively towards the lead car and eventually ends in the red region of Fig.~\ref{fig:lmpc-results-bridge} at the same time, resulting in significant packet loss and increased braking. 
With SR-LMPC, in early iterations the controller attempts to avoid this situation but saturates the control inputs. 
As it iterates, SR-LMPC begins to avoid saturation, or it saturates for a shorter amount of time. Overall, SR-LMPC is able to minimize fuel cost while also maintaining closer headway (top of Fig.~\ref{fig:lmpc-results-bridge} as well as bottom) and constant communication.

The intuition behind these results is as follows. Without the ability to predict communication performance, the following vehicle attempts to optimize on following distance and terminal cost. Once it enters the dead zone, it starts to drop packets and can only use increasingly shorter portions of (previously communicated) leader trajectory predictions. In the case of SR-LMPC with communication prediction, the algorithm is able to access a generic, unknown function that tells it whether the channel is expected to change in the future. SR-LMPC is able to converge on a trajectory that smoothly slows down in advance of the dead zone, which results in significantly improved global performance.

\section{Conclusions}\label{sec:conclusions}
In this paper, an extension to learning Model Predictive Control (LMPC) is presented. The controller is designed for applications to motion planning in dynamic environments, particularly when one or more of the decision variables comes from black box or data-driven models. The control architecture leverages a nominal outer loop motion planner, and then iterates over this trajectory candidate in an inner loop to find optimal policies in terms of both model- and data-driven variables. This outer and inner control scheme proceeds in a receding horizon fashion until the system reaches its objectives. These concepts are applied to connected autonomous vehicles and the notion of platooning, or collaborative adaptive cruise control.

To demonstrate the approach, a simulation of a leader-follower scenario for two connected autonomous vehicles is developed. The scenario includes physical characteristics that cause uncertainty in the communication channel, and the controller leverages recent advances in wireless channel prediction using machine learning. The SR-LMPC framework is able to generate improved trajectories in terms of not only communication, but also energy efficiency.

\addtolength{\textheight}{-12cm}   







\bibliographystyle{IEEEtran}
\bibliography{main.bib}

\begin{thebibliography}{10}
\providecommand{\url}[1]{#1}
\csname url@samestyle\endcsname
\providecommand{\newblock}{\relax}
\providecommand{\bibinfo}[2]{#2}
\providecommand{\BIBentrySTDinterwordspacing}{\spaceskip=0pt\relax}
\providecommand{\BIBentryALTinterwordstretchfactor}{4}
\providecommand{\BIBentryALTinterwordspacing}{\spaceskip=\fontdimen2\font plus
\BIBentryALTinterwordstretchfactor\fontdimen3\font minus
  \fontdimen4\font\relax}
\providecommand{\BIBforeignlanguage}[2]{{%
\expandafter\ifx\csname l@#1\endcsname\relax
\typeout{** WARNING: IEEEtran.bst: No hyphenation pattern has been}%
\typeout{** loaded for the language `#1'. Using the pattern for}%
\typeout{** the default language instead.}%
\else
\language=\csname l@#1\endcsname
\fi
#2}}
\providecommand{\BIBdecl}{\relax}
\BIBdecl

\bibitem{nowzari2016distributed}
C.~Nowzari and J.~Cort{\'e}s, ``Distributed event-triggered coordination for
  average consensus on weight-balanced digraphs,'' \emph{Automatica}, vol.~68,
  pp. 237--244, 2016.

\bibitem{dimarogonas2012distributed}
D.~V. Dimarogonas, E.~Frazzoli, and K.~H. Johansson, ``Distributed
  event-triggered control for multi-agent systems,'' \emph{IEEE Transactions on
  Automatic Control}, vol.~57, no.~5, pp. 1291--1297, 2012.

\bibitem{wang2016multi}
L.~Wang, A.~D. Ames, and M.~Egerstedt, ``Multi-objective compositions for
  collision-free connectivity maintenance in teams of mobile robots,'' in
  \emph{Decision and Control (CDC), 2016 IEEE 55th Conference on}.\hskip 1em
  plus 0.5em minus 0.4em\relax IEEE, 2016, pp. 2659--2664.

\bibitem{wang2016safety}
L.~Wang, A.~Ames, and M.~Egerstedt, ``Safety barrier certificates for
  heterogeneous multi-robot systems,'' in \emph{American Control Conference
  (ACC), 2016}.\hskip 1em plus 0.5em minus 0.4em\relax IEEE, 2016, pp.
  5213--5218.

\bibitem{luo2010model}
L.-h. Luo, H.~Liu, P.~Li, and H.~Wang, ``Model predictive control for adaptive
  cruise control with multi-objectives: comfort, fuel-economy, safety and
  car-following,'' \emph{Journal of Zhejiang University SCIENCE A}, vol.~11,
  no.~3, pp. 191--201, 2010.

\bibitem{corona2008adaptive}
D.~Corona and B.~De~Schutter, ``Adaptive cruise control for a smart car: A
  comparison benchmark for mpc-pwa control methods,'' \emph{IEEE Transactions
  on Control Systems Technology}, vol.~16, no.~2, pp. 365--372, 2008.

\bibitem{varaiya1993smart}
P.~Varaiya, ``Smart cars on smart roads: problems of control,'' \emph{IEEE
  Transactions on automatic control}, vol.~38, no.~2, pp. 195--207, 1993.

\bibitem{iihoshi2000vehicle}
A.~Iihoshi, S.~Kobayashi, and Y.~Furukawa, ``Vehicle platoon control system,''
  Feb.~29 2000, uS Patent 6,032,097.

\bibitem{hedrick1991longitudinal}
J.~Hedrick, D.~McMahon, V.~Narendran, and D.~Swaroop, ``Longitudinal vehicle
  controller design for ivhs systems,'' in \emph{American Control Conference,
  1991}.\hskip 1em plus 0.5em minus 0.4em\relax IEEE, 1991, pp. 3107--3112.

\bibitem{kato2002vehicle}
S.~Kato, S.~Tsugawa, K.~Tokuda, T.~Matsui, and H.~Fujii, ``Vehicle control
  algorithms for cooperative driving with automated vehicles and intervehicle
  communications,'' \emph{IEEE Transactions on Intelligent Transportation
  Systems}, vol.~3, no.~3, pp. 155--161, 2002.

\bibitem{schmied2015nonlinear}
R.~Schmied, H.~Waschl, R.~Quirynen, M.~Diehl, and L.~del Re, ``Nonlinear mpc
  for emission efficient cooperative adaptive cruise control,''
  \emph{IFAC-PapersOnLine}, vol.~48, no.~23, pp. 160--165, 2015.

\bibitem{stanger2013model}
T.~Stanger and L.~del Re, ``A model predictive cooperative adaptive cruise
  control approach,'' in \emph{2013 American Control Conference}.\hskip 1em
  plus 0.5em minus 0.4em\relax IEEE, 2013, pp. 1374--1379.

\bibitem{sancar2014mpc}
F.~E. Sancar, B.~Fidan, J.~P. Huissoon, and S.~L. Waslander, ``Mpc based
  collaborative adaptive cruise control with rear end collision avoidance,'' in
  \emph{2014 IEEE Intelligent Vehicles Symposium Proceedings}.\hskip 1em plus
  0.5em minus 0.4em\relax IEEE, 2014, pp. 516--521.

\bibitem{kim2006accurate}
K.-H. Kim and K.~G. Shin, ``On accurate measurement of link quality in
  multi-hop wireless mesh networks,'' in \emph{Proceedings of the 12th annual
  international conference on Mobile computing and networking}.\hskip 1em plus
  0.5em minus 0.4em\relax ACM, 2006, pp. 38--49.

\bibitem{de2005high}
D.~S. De~Couto, D.~Aguayo, J.~Bicket, and R.~Morris, ``A high-throughput path
  metric for multi-hop wireless routing,'' \emph{Wireless networks}, vol.~11,
  no.~4, pp. 419--434, 2005.

\bibitem{gnawali2009collection}
O.~Gnawali, R.~Fonseca, K.~Jamieson, D.~Moss, and P.~Levis, ``Collection tree
  protocol,'' in \emph{Proceedings of the 7th ACM conference on embedded
  networked sensor systems}.\hskip 1em plus 0.5em minus 0.4em\relax ACM, 2009,
  pp. 1--14.

\bibitem{farkas2006pattern}
K.~Farkas, T.~Hossmann, L.~Ruf, and B.~Plattner, ``Pattern matching based link
  quality prediction in wireless mobile ad hoc networks,'' in \emph{Proceedings
  of the 9th ACM international symposium on Modeling analysis and simulation of
  wireless and mobile systems}.\hskip 1em plus 0.5em minus 0.4em\relax ACM,
  2006, pp. 239--246.

\bibitem{liu2011foresee}
T.~Liu and A.~E. Cerpa, ``Foresee (4c): Wireless link prediction using link
  features,'' in \emph{Information Processing in Sensor Networks (IPSN), 2011
  10th International Conference on}.\hskip 1em plus 0.5em minus 0.4em\relax
  IEEE, 2011, pp. 294--305.

\bibitem{boano2010triangle}
C.~A. Boano, M.~Zuniga, T.~Voigt, A.~Willig, and K.~R{\"o}mer, ``The triangle
  metric: Fast link quality estimation for mobile wireless sensor networks,''
  in \emph{International Conference on Computer Communication Networks, 2010,
  Zurich, Switzerland}, 2010.

\bibitem{elnaggarbayesian}
M.~Elnaggar, K.~Whitehouse, and C.~H. Fleming, ``Bayesian wireless channel
  prediction for safety-critical connected autonomous vehicles.''

\bibitem{rosolia2018learning}
U.~Rosolia and F.~Borrelli, ``Learning model predictive control for iterative
  tasks. a data-driven control framework,'' \emph{IEEE Transactions on
  Automatic Control}, vol.~63, no.~7, pp. 1883--1896, 2018.

\bibitem{brunner2017repetitive}
M.~Brunner, U.~Rosolia, J.~Gonzales, and F.~Borrelli, ``Repetitive learning
  model predictive control: An autonomous racing example,'' in \emph{2017 IEEE
  56th Annual Conference on Decision and Control (CDC)}.\hskip 1em plus 0.5em
  minus 0.4em\relax IEEE, 2017, pp. 2545--2550.

\bibitem{rosolia2017robust}
U.~Rosolia, X.~Zhang, and F.~Borrelli, ``Robust learning model predictive
  control for iterative tasks: Learning from experience,'' in \emph{2017 IEEE
  56th Annual Conference on Decision and Control (CDC)}.\hskip 1em plus 0.5em
  minus 0.4em\relax IEEE, 2017, pp. 1157--1162.

\bibitem{firoozi2018safe}
R.~Firoozi, S.~Nazari, J.~Guanetti, R.~O'Gorman, and F.~Borrelli, ``Safe
  adaptive cruise control with road grade preview and v2v communication,''
  \emph{arXiv preprint arXiv:1810.09000}, 2018.

\bibitem{rosolia2017learning}
U.~Rosolia and F.~Borrelli, ``Learning model predictive control for iterative
  tasks: A computationally efficient approach for linear system,''
  \emph{IFAC-PapersOnLine}, vol.~50, no.~1, pp. 3142--3147, 2017.

\bibitem{jafarzadeh2018exact}
H.~Jafarzadeh and C.~H. Fleming, ``An exact geometry--based algorithm for path
  planning,'' \emph{International Journal of Applied Mathematics and Computer
  Science}, vol.~28, no.~3, pp. 493--504, 2018.

\bibitem{ye1989extension}
Y.~Ye and E.~Tse, ``An extension of karmarkar's projective algorithm for convex
  quadratic programming,'' \emph{Mathematical programming}, vol.~44, no. 1-3,
  pp. 157--179, 1989.

\bibitem{bashiri2017platoon}
M.~Bashiri and C.~H. Fleming, ``A platoon-based intersection management system
  for autonomous vehicles,'' in \emph{2017 IEEE Intelligent Vehicles Symposium
  (IV)}.\hskip 1em plus 0.5em minus 0.4em\relax IEEE, 2017, pp. 667--672.

\bibitem{bashiri2018paim}
M.~Bashiri, H.~Jafarzadeh, and C.~H. Fleming, ``Paim: Platoon-based autonomous
  intersection management,'' in \emph{2018 21st International Conference on
  Intelligent Transportation Systems (ITSC)}.\hskip 1em plus 0.5em minus
  0.4em\relax IEEE, 2018, pp. 374--380.

\bibitem{schneider2016directional}
C.~Schneider, M.~K{\"a}ske, G.~Sommerkorn, R.~S. Thom{\"a}, A.~Roivainen,
  J.~Meinil{\"a}, and V.~Tervo, ``Directional analysis of multipath propagation
  in vehicle-2-vehicle channels,'' in \emph{2016 10th European Conference on
  Antennas and Propagation (EuCAP)}.\hskip 1em plus 0.5em minus 0.4em\relax
  IEEE, 2016, pp. 1--5.

\bibitem{cheng2013roadside}
L.~Cheng, D.~D. Stancil, and F.~Bai, ``A roadside scattering model for the
  vehicle-to-vehicle communication channel,'' \emph{IEEE Journal on Selected
  Areas in Communications}, vol.~31, no.~9, pp. 449--459, 2013.

\bibitem{Hs2006}
D.~O.~T. Hs, ``{Vehicle Safety Communications Project - Final Report},''
  \emph{Communications}, no. April, p.~44, 2006.

\bibitem{cplex2009v12}
I.~I. CPLEX, ``V12. 1: User’s manual for cplex,'' \emph{International
  Business Machines Corporation}, vol.~46, no.~53, p. 157, 2009.

\end{thebibliography}

\end{document}